\documentclass[12pt]{amsart}
\usepackage{amssymb}

\newtheorem{thm}{Theorem}[section]

\newtheorem{lem}[thm]{Lemma}

\newtheorem{prop}[thm]{Proposition}

  \def\Cset{\mathbb{C}}
  
  \def\Nset{\mathbb{N}}
  
  \def\Rset{\mathbb{R}}

\newcommand{\subsetnum}[2]{\genfrac{\{}{\}}{0pt}{}{#1}{#2}}

\begin{document}

\title{Contributions to the Theory of the Barnes Function}

\author{V. S. Adamchik}

\address{Carnegie Mellon University \\Pittsburgh, USA}

\email{adamchik@cs.cmu.edu}

\thanks{This work was supported by grant CCR-0204003 from the National
Science Foundation.}


\begin{abstract}
This paper presents a family of new integral representations and
asymptotic series of the multiple gamma function. The numerical
schemes for high-precision computation of the Barnes gamma
function and Glaisher's constant are also discussed.
\end{abstract}


\maketitle
\section{Introduction}

In a sequence of papers published between 1899-1904, Barnes
introduced and studied (see
\cite{Barnes99,Barnes00,Barnes01,Barnes04}) a generalization of
the classical Euler gamma function, called the multiple gamma
function ${\Gamma}_{n}(z)$. The function \({{\Gamma }_n}(z)\)
satisfies the following recurrence-functional equation
\cite{Vardi88,Vign}:

\begin{equation}
\label{GDef}
\begin{array}{l}
\displaystyle
{\Gamma}_{n+1}(z+1) =
\frac{{\Gamma}_{n+1}(z)}{{\Gamma}_{n}(z)}, \quad z \in
\Cset, \quad n \in \Nset,  \\
\noalign{\vskip 1.0pc} \displaystyle
{\Gamma}_{1}(z) = {\Gamma}(z), \\
\noalign{\vskip 1.0pc} \displaystyle
 {\Gamma}_{n}(1) = 1,
\end{array}
\end{equation}
\noindent
where ${\Gamma}(z)$ is the Euler gamma function. The system
(\ref{GDef}) has a unique solution if we add the condition of
convexity \cite{Vign}

\[
(-1)^{n+1} \frac{d^{n+1}}{d \, z^{n+1}}\, \log \Gamma_{n}(z) \geq 0,
\quad z > 0.
\]

In this paper we aim at the most interesting case $G(z) =
1/\Gamma_{2}(z)$, being the so-called double gamma function or the
Barnes $G$-function:
\begin{equation}
\begin{array}{l}
\displaystyle
G(z+1) = \Gamma(z) \> G(z), \quad z \in \Cset, \\
\noalign{\vskip 1.0pc} \displaystyle
G(1) = G(2) = G(3) = 1, \nonumber \\
\noalign{\vskip 1.0pc} \displaystyle
 \frac{d^{3}}{d \, z^{3}}\, \log G(z) \geq 0, \quad z
> 0.
\end{array}
\end{equation}

The $G$-function has several equivalent forms
including the Weierstrass canonical product
\[
G(z+1)={{(2 \pi )}^{\frac{z}{2} }} \exp\Big(-\frac{z+{z^2} (1+\gamma)}{2}
\Big) \prod _{k=1}^{\infty }{{\Big(1+\frac{z}{k}\Big)}^k}
\exp\Big(\frac{{z^2}}{2 k}-z\Big)
\]
and the functional relationship with the Hurwitz zeta function
$\zeta (t,z)$, due to Gosper \cite{Gosper} and Vardi \cite{Vardi88}

\begin{equation}
\label{G2Z} \log G(z+1)-z \, \log  \Gamma (z)={{\zeta }^{\prime
}}(-1)-{{\zeta }^{\prime }}(-1,z), \quad \Re(z)>0,
\end{equation}

\smallskip
\noindent
where \(\gamma \) denotes the Euler-Mascheroni constant and the derivative
${\zeta }^{\prime }(t,z) = \frac{d}{dt} \zeta(t, z)$ is
defined as an analytic continuation provided by the Hermite integral (\ref{hermite}):

\begin{equation}
\label{zeta(-1, z)}
\begin{array}{l}
\displaystyle
{{\zeta }^{\prime}}(-1,z) = \frac{z^2}{2} \log z -
\frac{z^2}{4} - \frac{z}{2} \log z  \> + \\
\noalign{\vskip 1.0pc} \displaystyle
 2 \, z \int _{0}^{\infty}
\frac{\arctan (\frac{x}{z})} {e^{2 \pi x} - 1} \> dx +
  \int_{0}^{\infty} \frac{x\,\log (x^2 + z^2)}{e^{2 \pi x} - 1} \> dx.
  \quad \Re(z)>0,
\end{array}
\end{equation}

\medskip
\noindent
In particular,

\begin{equation}
\label{zeta(-1)}
{{\zeta }^{\prime }}(-1) = {{\zeta
}^{\prime}}(-1,0) = 2 \int_{0}^{\infty}
\frac{x\,\log x}{e^{2 \pi x} - 1} \> dx.
\end{equation}
\noindent For \(z\rightarrow  \infty \), the $G$-function
has the Stirling asymptotic expansion

\[
\log  G(z+1)  = \frac{z}{2} \log(2 \pi) - \log A -
\frac{3{z^2}}{4} + \bigg(\frac{{z^2}}{2}-\frac{1}{12}\bigg) \log z
+ O\Big(\frac{1}{z}\Big),
\]
where $A$ is the Glaisher-Kinkelin constant given by
\begin{equation}
\label{GK}
\log A = \frac{1}{12}-\zeta^{\prime}(-1).
\end{equation}

\noindent Combining this definition with the representation
(\ref{zeta(-1)}), we obtain a new integral for $\log A$:

\begin{equation}
\log A = \frac{1+\log(2 \pi)}{12} - \frac{1}{2 \pi^2}
\int_{0}^{\infty} \frac{x \log x}{e^{x}-1} \> dx.
\end{equation}

\noindent The constant $A$ originally appeared in papers by
Kinkelin \cite{Kinkelin} and Glaisher
\cite{Glaisher77,Glaisher87,Glaisher93} on the asymptotic
expansion (when $n \rightarrow \infty$) of the following product

\[
1^{1^{p}} \, 2^{2^{p}} \, \ldots n^{n^{p}}, \quad p \in \Nset.
\]

\noindent Similar to the Euler gamma function, the $G$-function satisfies the multiplication formula \cite{Vardi88}:

\begin{equation}
\label{multipl}
\begin{array}{c}
\displaystyle
G(n \> z)=e^{{{\zeta }^{\prime}}(-1) (1-{n^2})} n^{n^2 z^2/2-n z +5/12} \\
\noalign{\vskip 0.5pc} \displaystyle
 {{(2 \pi )}^{\frac{n-1}{2}(1-n z)}} \prod _{i=0}^{n-1} \, \prod
_{j=0}^{n-1}G\Big(z+\frac{i+j}{n}\Big), \quad n \in \Nset \\
\end{array}
\end{equation}

    The Barnes function is not listed in the tables of the
most well-known special functions. However, it is cited in the
exercises by Whittaker and Watson \cite{WW}, and in entries 6.441,
8.333 by Gradshteyn and Ryzhik \cite{GR}. Recently, Vardi
\cite{Vardi88}, Quine and Choi \cite{Quine}, and Kumagai
\cite{Kumagai} computed the functional determinants of Laplacians
of the $n$-sphere in terms of the \({{\Gamma }_n}\)-function.
Choi, Srivastava \cite{Choi95,Choi99} and Adamchik
\cite{ISSAC01,AMC03} presented relationships (including integrals
and series) between other mathematical functions and the multiple
gamma function.

\medskip
    The \({{\Gamma }_n}\) function has several interesting applications in pure and applied mathematics and theoretical physics. The most intriguing one is the application of $G$-function to the Riemann Hypothesis. Montgomery \cite{Montgomery} and  Sarnak \cite{Sarnak97} (see \cite{Keating00} for additional references) have conjectured that the non-trivial zeros of the Riemann zeta function are pairwise distributed like eigenvalues of matrices in the Gaussian Unitary Ensemble (GUE) of random matrix theory. It has been shown in works by Mehta, Sarnak, Conrey, Keating, and Snaith that a closed representation for statistical averages over GUE of $N \times N$ unitary matrices, when $N\rightarrow  \infty$ can be expressed in terms of the Barnes functions.

\medskip
    Another interesting appearance of the $G$-function is in its relation to the determinants of the Hankel matrices. Consider a matrix $M$ of the Bell numbers:

\begin{displaymath}
 {M_n}=\left(\begin{array}[c]{cccc}
    {B_1}&{B_{2}}&...&{B_{n}} \\
    {B_2}&{B_{3}}&...&{B_{n+1}} \\
      &  &...&   \\
    {B_n}&{B_{n+1}}&...&{B_{2n-1}}
  \end{array}\right)
\end{displaymath}
where $B_{n}$ are the Bell numbers, defined by

\[
B_n=\sum_{k=0}^{n} {\subsetnum{n}{k}}
\]

\noindent
and $\subsetnum{n}{k}$ are the Stirling subset numbers \cite{Rosen}

\begin{equation}
\label{subset} {\subsetnum{n}{k}} = k \, {\subsetnum{n-1}{k}} +
{\subsetnum{n-1}{k-1}} , \quad {\subsetnum{n}{0}} =
{\bigg\{\begin{array}[c]{c}
    1, \quad n = 0, \\ \noalign{\vskip 0.25pc}
    0, \quad n \neq 0
  \end{array}}
\end{equation}

\smallskip
\noindent
It was proved in \cite{Aigner}, that

\[ \det M_n =G(n+1),\quad n \in \Nset.\]

\bigskip

\section{Special cases of the Barnes Function}

For some particular arguments, the Barnes function can be
expressed in terms of the known special functions and constants. The simplest
(but not-trivial) special case is due to Barnes \cite{Barnes99}:

\begin{equation}
\label{G2(1/2)} \log {\rm
G}\Big(\frac{1}{2}\Big)=\frac{1}{8}+\frac{\log 2}{24}-\frac{\log
\pi }{4}-\frac{3 }{2} \log A,
\end{equation}

\smallskip
\noindent
where $A$ is the Glaisher-Kinkelin constant (\ref{GK}). The next
two special cases are due to Srivastava and Choi \cite{Choi99}:

\begin{equation}
\label{G(1/4)} \log G\Big(\frac{1}{4}\Big) =
\frac{3}{32}-\frac{{\bf G}}{4 \pi }-\frac{3}{4} \log  \Gamma
\Big(\frac{1}{4}\Big) - \frac{9}{8} \log  A,
\end{equation}

\begin{equation}
\label{G(3/4)} \log G\Big(\frac{3}{4}\Big) = \log {\rm
G}\Big(\frac{1}{4}\Big)+\frac{{\bf G}}{2  \pi }-\frac{\log
2}{8}-\frac{\log  \pi }{4}+ \log \Gamma \Big(\frac{1}{4}\Big),
\end{equation}

\smallskip
\noindent
where ${\bf G}$ denotes Catalan's constant:

\[
{\bf G} = \sum_{k=0}^{\infty} \frac{(-1)^k}{(2k+1)^2}.
\]

\smallskip
\noindent
Identities (\ref{G(1/4)}) and (\ref{G(3/4)}) are proven by using
the duplication formula (\ref{multipl}) (setting $n=2$) together
with (\ref{G2(1/2)}). More special cases for $z={\frac{1}{6}},
{\frac{1}{3}},{\frac{2}{3}}$ and ${\frac{5}{6}}$ can be derived as well.
Here is one of such formulas (see \cite{ISSAC01} for the others):

\begin{prop}
\begin{equation}
\label{log_g(1/3)}
\begin{split}
\log G\left(\frac{1}{3}\right)=&\frac{1}{9}+\frac{\log 3}{72}+\frac{\pi
}{18  {\sqrt{3}}}\> -\\
&\frac{2}{3} \log \Gamma\left(\frac{1}{3}\right)-\frac{4}{3}\log A-\frac{{{\psi
}^{(1)}}\big(\frac{1}{3}\big)}{12 {\sqrt{3}}  \pi },
\end{split}
\end{equation}

\noindent
where $\psi^{(1)}(z)=\frac{{\partial}^2}{\partial {z^2}}  \log
\Gamma (z)$ is the polygamma function.
\end{prop}
\smallskip
\begin{proof}[Proof.]
We recall the Lerch functional equation for the Hurwitz zeta
function (see \cite{JCAM98,Bateman})

\[
\begin{array}[c]{c}
\displaystyle
{\zeta}^{\prime}(-n,z) + (-1)^n {\zeta}^{\prime}(-n,1-z)  =
\frac{\pi i}{n+1} \> B_{n+1}(z) + \\
\noalign{\vskip0.5pc} \displaystyle {e^{-\pi i n/2}}
\frac{n!}{(2\pi)^n}\> {{\rm Li}_{n+1}} \big(e^{2 \pi i z}\big),
\quad 0<z<1,  \quad n \in \Nset,
\end{array}
\]

\smallskip
\noindent
where ${\zeta}^{\prime}(t,z) = \frac{d}{d t}\zeta (t,z)$,
$B_{n}(z)$ are the Bernoulli polynomials and ${\rm Li}_{n}(z)$ is
the polylogarithm. Setting $n=1$ and using (\ref{G2Z}) yields the
reflexion formula for the $G$-function (see also \cite{Choi99,AMC03}):

\begin{equation}
\label{reflex}
\log \Big(\frac{G(1+z)}{G(1-z)}\Big)=-z  \log \Big(\frac{\sin(\pi
z)}{\pi }\Big)-\frac{1}{2  \pi }{{\rm Cl}_{2}}(2 \,\pi z), \quad 0<z<1,
\end{equation}

\smallskip
\noindent
where ${\rm Cl}_{2}(z)$ is the Clausen function defined by

\begin{equation}
\label{clausen}
{\rm Cl}_{2}(z) = -\Im ({\rm Li}_{2}(e^{-i z})).
\end{equation}

\smallskip
\noindent
Upon setting $z=\frac{1}{3}$ to (\ref{reflex}), and making use of

\[
{\rm Cl}_{2}\Big(\frac{2\pi }{3}\Big)=\frac{{{\psi
}^{(1)}}\big(\frac{1}{3}\big)}{3  {\sqrt{3}}}-\frac{2  {{\pi
}^2}}{9  {\sqrt{3}}},
\]
we obtain

\begin{equation}
\label{g(1/3)}
\frac{G\left(\frac{1}{3}\right)}{G\left(\frac{2}{3}\right)}=\frac{
{\root{3}\of{2  \pi }}}{{\root{6}\of{3}} \, \Gamma
\big(\frac{1}{3}\big)}\, {\exp\bigg({\frac{2  {{\pi }^2}-3  {{\psi
}^{(1)}}\left(\frac{1}{3}\right)}{18  \pi   {\sqrt{3}}  }}\bigg)}.
\end{equation}

\noindent
On the other hand, using the multiplication formula (\ref{multipl}) with \(z=\frac{1}{3}\)
and $n=3$, we find

\begin{equation}
\label{g(1/3)2} G\Big(\frac{1}{3}\Big) {\rm
G}\Big(\frac{2}{3}\Big)={\root{3}\of{\frac{{3^{7/12}} {e^{2/3}}}{2
\pi   {A^8}  \Gamma \big(\frac{1}{3}\big)}}}.
\end{equation}
Now, combining (\ref{g(1/3)}) and (\ref{g(1/3)2}), we conclude the
proof.
\end{proof}

It remains to be seen if the Barnes function (or the multiple
gamma function) can be expressed as a finite combination of
elementary functions and constants for other rational arguments.
The most general result available for the multiple gamma function
$\Gamma_{n}(z)$ is a closed form at $z = \frac{1}{2}$:

\begin{prop}
For $n \in \Nset$
\begin{equation}
\label{Gamma.half}
\begin{array}{l}
\displaystyle

(-1)^n (n-1)! \log \Gamma_{n}(\frac{1}{2}) = - \frac{(2n-3)!! \log
\pi}{2^n} \> + \\
\noalign{\vskip 1.0pc}{\hskip 3.0pc} \displaystyle

\log 2 \sum_{k=1}^{n}\frac{P_{k,n}(\frac{1}{2}) \,
B_{k+1}}{(k+1)\, 2^{k}} +

\sum_{k=1}^{n} \frac{2^{k}-1}{2^k} P_{k,n}(\frac{1}{2})
\zeta^{\prime}(-k),
\end{array}
\end{equation}
where $P_{k,n}(\frac{1}{2})$ are coefficients by $x^{k}$, $k\leq
n$ in expansion of

\[
 \prod_{j=1}^{n-1} (x + j - \frac{1}{2})
\]
\end{prop}

\begin{proof}
The proof follows directly from \cite{lanl03}, formula (17), by
setting $z=\frac{1}{2}$.
\end{proof}

Here are a few particular cases:

the gamma function ($n = 1$):

\[
\Gamma(\frac{1}{2}) = \sqrt{\pi}
\]

\medskip
the Barnes function ($n = 2$):

\[
G(\frac{1}{2}) = \frac{2^{1/24} \, e^{1/8}}{A^{3/2} \, \pi^{1/4}}
\]

\medskip
the triple gamma function ($n = 3$):

\[
\Gamma_{3}(\frac{1}{2}) = \frac{A^{3/2} \, \pi^{3/16}}{2^{1/24}}
\, \exp\left( \frac{7 \, \zeta(3)}{32 \, \pi^2} -
\frac{1}{8}\right)
\]
\bigskip

\section{Hermite integrals}

In this section we derive a few integral representations for the
Barnes function. We start with recalling the Hermite integral (see
\cite{AAR,Magnus}) for the Hurwitz zeta function:

\begin{equation}
\label{hermite} \zeta(s,z) =
\frac{{z^{-s}}}{2}+\frac{{z^{1-s}}}{s-1}+2\int _{0}^{\infty
}\frac{\sin(s  \arctan(x/z))}{{{({x^2}+{z^2}) }^{s/2}}({e^{2 \pi
x}}-1)} \>  dx,
\end{equation}
which provides analytic continuation of $\zeta(s,z)$ to the domain
$s \in \Cset-\{1\}$. Differentiating both sides of (\ref{hermite})
with respect to $s$, letting $s=-1$ (and \(s=-2\)), and using the
second Binet formula for $\log \Gamma(z)$ (see \cite{AAR})

\begin{equation}
\label{logGamma} \log \Gamma(z) = (z- \frac{1}{2}) \log z - z +
\frac{\log (2 \pi)} {2}  +2\int _{0}^{\infty
}\frac{\arctan(x/z)}{{e^{2 \pi x}}-1} \> dx
\end{equation}
together with (\ref{G2Z}), we readily obtain:

\begin{equation}
\label{hermite_for_G}
\begin{array}{c}
\displaystyle
\log G(z+1) = \frac{{z^2}}{2}\Big(\log z-{H_2}\Big) -\Big(z \, {{\zeta
}^{\prime }}(0)-  {{\zeta }^{\prime }}(-1)\Big) \, - \\
\noalign{\vskip 1.0pc} \displaystyle
\int_{0}^{\infty }\frac{x
\log({x^2}+{z^2})}{{e^{2  \pi   x}}-1} \>  dx,  \quad  \Re(z)>0,
\end{array}
\end{equation}

\smallskip
\noindent
where \({H_p}=\sum _{k=1}^{p}{1/k}\) are the harmonic numbers and  \({{\zeta }^{\prime }}(z)\) is a derivative of the Riemann zeta function.

\medskip
In the same manner, repeating the above steps $n$ times,  we
derive the general integral representations.

\medskip
\begin{prop}
For \(\Re(z)>0\)

\begin{equation}
\label{prop6.1}
\begin{array}{l}
\displaystyle
(2 n)! \, \log  {{\Gamma}_{2n+1}}(z+1) =        \\
\noalign{\vskip 1.0pc} \displaystyle {\hskip 2.5pc}
\frac{z^{2 n+1} (\log z -
H_{2 n+1})}{2 n+1} - \sum_{k=0}^{2 n}(-1)^k \binom{2n}{k} \,
{\zeta^{\prime}}(-k) \> {z^{2 n-k}} \>-  \\
  \noalign{\vskip 1.0pc} \displaystyle {\hskip 2.5pc}
\sum_{k=1}^{2 n-1}{(-1)^k}  k! \> {\subsetnum{2 n}{k}} \log
{{\Gamma}_{k+1}}(z+1) \>+  \\
  \noalign{\vskip 1.0pc} \displaystyle {\hskip 2.5pc}
2 \> {(-1)^n}  \int_{0}^{\infty }
\frac{x^{2 n} \arctan\big({x/z}\big)}{e^{2 \pi x}-1}\, dx
\end{array}
\end{equation}

\medskip
\noindent and

\begin{equation}
\label{prop6.2}
\begin{array}{l}
\displaystyle
(2 n-1)!  \log  {{\Gamma }_{2n}}(z+1) = -\frac{{z^{2  n}}  (\log  z-{H_{2  n}})}{2 n} \>+         \\
  \noalign{\vskip 1.0pc} \displaystyle {\hskip 2.5pc}

\sum _{k=0}^{2
n-1}{{(-1)}^k}  \binom{2n-1}{k} {{\zeta }^{\prime }}(-k)  {z^{2n-k-1}}\>+ \\
\noalign{\vskip 1.0pc} \displaystyle {\hskip 2.5pc}

\sum _{k=1}^{2n-2}{(-1)^k} k!\, {\subsetnum{2 n-1}{k}} \log  {{\Gamma }_{k+1}}(z+1)\>-\\
\noalign{\vskip 1.0pc} \displaystyle {\hskip 2.5pc}

{{(-1)}^n}\int _{0}^{\infty }\frac{{x^{2  n-1}}
  \log({x^2}+{z^2})}{{e^{2  \pi   x}}-1}\,dx\\\noalign{\vskip 0.5pc}
\end{array}
\end{equation}
where $\subsetnum{n}{k}$ are the Stirling subset numbers, defined
in (\ref{subset}).
\end{prop}

\medskip
Upon substituting $z=1$ into formulas (\ref{prop6.1})
and (\ref{prop6.2}), we obtain new integrals

\begin{equation}
\label{int1} 2\, {(-1)^n}  \int_{0}^{\infty
}\frac{{x^{2n}}\arctan(x)}{{e^{2  \pi x}}-1}\,
dx=\frac{{H_{2n+1}}}{2n+1}+\sum_{k=0}^{2n}{{(-1)}^k}
\binom{2n}{k} {{\zeta }^{\prime }}(-k)
\end{equation}

\begin{equation}
\label{int2} {(-1)}^n  \int_{0}^{\infty}\frac{{x^{2n-1}}
\log(1+{x^2})}{{e^{2  \pi x}}-1}\,dx = \frac{{H_{2n}}}{2n}+
\sum_{k=0}^{2n-1}{(-1)^k \binom{2n-1}{k} {{\zeta }^{\prime }}(-k)}
\end{equation}

Performing simple transformations of the integrands in (\ref{int1}) and (\ref{int2}) will apparently lead to a large number of special integrals:

\[ 2\int _{0}^{\infty }\frac{x \log x}{{e^{2  \pi   x}}-1}\, dx =  {{\zeta }^{\prime }}(-1) \]

\[ \int _{0}^{\infty }\frac{x \log
x}{{e^{2  \pi   x}}+1}\, dx =  12  {{\zeta }^{\prime }}(-1)  + \log
2 \]

\[ \int _{0}^{\infty }\frac{x
\log(1+{x^2})}{{e^{2  \pi   x}}+1}\, dx = \frac{3}{4}-\frac{23}{24}
  \log  2+\frac{1}{2}  {{\zeta
}^{\prime }}(-1) \]

\[ 2  \int _{0}^{\infty
}\frac{\arctan\big(\frac{x}{z}\big)}{{e^{2  \pi x}}+1}\, dx = z
\log z-z+\frac{\log(2  \pi )}{2}  -\log \Gamma
\big(z+\frac{1}{2}\big), \quad \Re(z)>0 \]

\[ 4  \int _{0}^{\infty}\frac{\arctan(x)}{{e^{2  \pi   x}}+1}\, dx = -2 + 3 \log 2 \]

\[2  \int _{0}^{\infty }\frac{x \, dx}{({e^{2  \pi   x}}+1)  ({x^2}+{z^2})}=\psi
\big(z+\frac{1}{2}\big)-\log  z, \quad
\arg(z)\neq \frac{\pi }{2} \]

\[2  \int _{0}^{\infty}\frac{x \, dx}{({e^{2  \pi   x}}-1) ({x^2}+{z^2})}=\log  z-\psi
(z)-\frac{1}{2 z}, \quad
    \arg(z)\neq \frac{\pi}{2}\]

\bigskip
\bigskip
\section{Binet-like representation}

In this section we derive the Binet integral representation for
the Barnes function.

\begin{prop}
The Barnes G-function admits the Binet integral representation:

\begin{equation}
\label{binet}
\begin{array}{l}
\displaystyle
\log G(z+1) = z \log  \Gamma (z) + \frac{{z^2}}{4}-\frac{\log
z}{2} \, {B_2}(z) - \log A \, - \\
\noalign{\vskip0.5pc}{\hskip2.5pc} \displaystyle
\int_{0}^{\infty}\frac{{e^{-z  x}}}{{x^2}}
\Big(\frac{1}{1-{e^{-x}}}-\frac{1}{x}-\frac{1}{2}-\frac{x}{12}\Big)\, dx,\quad
                \Re(z)>0.
\end{array}
\end{equation}
\end{prop}

\medskip
\begin{proof}[Proof.]
Recall the well-known integral for the Hurwitz function
\cite{Bateman}

\begin{equation}
\label{hurwitz_int} \zeta (s,z)=\frac{1}{\Gamma
(s)}\int_{0}^{\infty }\frac{{x^{ s-1}} {e^{-z x}}}{1-{e^{-x}}} \>
dx, \quad \Re(s)>1,  \quad  \Re(z)>0
\end{equation}
The integral in the right hand-side of (\ref{hurwitz_int}) can be
analytically continued to the larger domain $\Re(s) > -2$ by
subtracting the truncated Taylor series of $\frac{1}{1-e^{-x}}$ at
$x=0$. Since

\begin{equation}
\label{taylor}
\frac{1}{1-e^{-x}} = \frac{1}{x}+\frac{1}{2}+\frac{x}{12} + O ( x^{3})
\end{equation}
and

\[
\begin{array}{c}
\displaystyle
\frac{1}{\Gamma(s)} \int_{0}^{\infty}{x^{s-1} e^{-z x}
\Big(\frac{1}{x}+\frac{1}{2}+\frac{x}{12} \Big)\, dx} =
\frac{s}{12} z^{-1-s} + \frac{z^{-s}}{2} - \frac{z^{1-s}}{1-s}, \\
\noalign{\vskip 0.5pc} \displaystyle
 \Re(s)>1, \Re(z)>0
\end{array}
\]
we therefore can rewrite (\ref{hurwitz_int}) as

\[
\begin{array}{l}
\displaystyle
\zeta (s,z)=-\frac{{z^{1-s}}}{1-s} + \frac{{z^{-s}}}{2} + \frac{s}{12} \,
{z^{-s-1}} \>+ \\
\noalign{\vskip1.0pc}{\hskip1.5pc} \displaystyle
\frac{1}{\Gamma (s)}\int _{0}^{\infty }{x^{s-1}} {e^{-z
x}}
\Big(\frac{1}{1-{e^{-x}}}-\frac{1}{x}-\frac{1}{2}-\frac{x}{12}\Big)\> dx, \quad \Re(s)>-2.
\end{array}
\]
Differentiating the above formula with respect to $s$ and computing the limit
at $s=-1$, we obtain

\[
\begin{array}{l}
\displaystyle
\zeta^{\prime}(-1,z)=\frac{1}{12}-\frac{{z^2}}{4}+\frac{\log
z}{2}\, {B_2}(z)+ \\
\noalign{\vskip 1.0pc}{\hskip2.5pc} \displaystyle
 \int _{0}^{\infty
}\frac{{e^{-z x}}}{{x^2}}
\Big(\frac{1}{1-{e^{-x}}}-\frac{1}{x}-\frac{1}{2}-\frac{x}{12}\Big)\,
dx, \quad \Re(z)>0,
\end{array}
\]
where $B_{2}(z) $ is the second Bernoulli polynomial. The
convergence of the above integral is ensured by (\ref{taylor}).
Finally, using (\ref{G2Z}), we conclude the proof.
\end{proof}

\bigskip

\section{Asymptotic expansion of  $\Gamma_{n}(z)$}

The asymptotic expansion for $\Gamma_{n}(z)$ when $z \rightarrow \infty$ follows
straightforwardly from formulas (\ref{prop6.1}) and (\ref{prop6.2}) by
expanding $\log(1+x^{2})$ and $\arctan(x)$ into the Taylor series, and then
performing formal term-by-term integration.

\begin{lem}
For $\Re(z)>0$

\begin{equation}
\label{lem.asym1}
\begin{array}{l}
\displaystyle
\log G(z+1) = \frac{z^2}{2}\left( \log z - \frac{3}{2}\right) - \frac{\log z}{12} - z \, \zeta^{\prime}(0) + \zeta^{\prime}(-1) \>- \\
\noalign{\vskip1.0pc}{\hskip6.0pc} \displaystyle

\sum_{k=1}^{n} \frac{B_{2k+2}}{4\,k\,(k+1)\,z^{2k}} +
O(\frac{1}{z^{2n+2}}),\quad z \rightarrow \infty.
\end{array}
\end{equation}
\end{lem}

\begin{proof}We start with the integral representation (\ref{hermite_for_G}).
By recalling the Taylor series for $\log(1+x)$

\[
\log(1 +\frac{x}{z}) = \sum_{k=0}^{n}\frac{(-1)^k}{k+1}\,
{(\frac{x}{z})}^{k+1} - \frac{(-1)^n}{z^{n+1}}
\int_{0}^{x}\frac{y^{n+1}}{y+z}\> dy,
\]

\noindent we find that

\[
\begin{array}{l}
\displaystyle

\int_{0}^{\infty} \frac{x \log(1+x^{2}/z^{2})}{e^{2\,
\pi x} - 1}\> dx = \\
\noalign{\vskip 1.0pc} \displaystyle

\sum_{k=0}^{n} \frac{(-1)^k}{k+1} \frac{1}{z^{2k+2}}
\int_{0}^{\infty} \frac{x^{2k+3}}{e^{2 \,\pi x} - 1}\> dx -
\frac{(-1)^n}{z^{2n+2}} \int_{0}^{\infty} \frac{t \>dt}{e^{2 \,\pi
t} - 1} \int_{0}^{t^{2}}\frac{ y^{n+1}}{y+z^{2}} \>dy
\end{array}
\]

\medskip \noindent which, by appealing to (formula 1.6.4 in
\cite{AAR})

\begin{equation}
\label{int.bern}
\int_{0}^{\infty} \frac{t^{2k-1} }{e^{2 \,\pi
\,t} - 1} \> dt = (-1)^{k+1} \> \frac{B_{2k}}{4 \, k}, \quad k \in
\Nset,
\end{equation}

\smallskip
\noindent reduces to

\[
\begin{array}{l}
\displaystyle

\int_{0}^{\infty} \frac{x \log(1+x^{2}/z^{2})}{e^{2\,
\pi x} - 1}\> dx = \\
\noalign{\vskip 1.0pc} \displaystyle

-\sum_{k=1}^{n+1} \frac{B_{2k+2}}{4 \,k \,(k+1) z^{2k}}  - 2\,
(-1)^n \int_{0}^{\infty} \frac{t \>dt}{e^{2 \,\pi \,t} - 1}
\int_{0}^{t/z}\frac{ y^{2n+3}}{1+y^{2}} \>dy
\end{array}
\]

\medskip \noindent Substituting this into (\ref{hermite_for_G}),
we complete the proof.
\end{proof}
\bigskip

\begin{lem}
For $\Re(z)>0$

\begin{equation}
\label{lem.asym2}
\begin{array}{l}
\displaystyle 2 \log \Gamma_{3}(z+1) = -\frac{z^3}{3}\Big( \log z
- \frac{11}{6}\Big) + \frac{z^2}{2}\Big(\log z - \frac{3}{2} +
\zeta^{\prime}(0)\Big) - \\
\noalign{\vskip1.0pc}{\hskip4.0pc} \displaystyle

\frac{\log z}{12} - z \Big( \zeta^{\prime}(0) + 2
\zeta^{\prime}(-1)\Big)  + \zeta^{\prime}(-1)+\zeta^{\prime}(-2) \>-\\
\noalign{\vskip1.0pc}{\hskip4.0pc} \displaystyle

\sum_{k=1}^{n} \frac{B_{2k+2}}{4\,k\,(k+1)\,z^{2k}} -
\sum_{k=1}^{n} \frac{B_{2k+2}}{2\,(k+1)\,(2k-1)\,z^{2k-1}} \> +\\
\noalign{\vskip1.0pc}{\hskip4.0pc} \displaystyle

O(\frac{1}{z^{2n+1}}),\quad z \rightarrow \infty.
\end{array}
\end{equation}
\end{lem}

\begin{proof}
Employing the same technique as in getting (\ref{lem.asym1}), we
begin with the integral representation of the triple gamma
function:

\begin{equation}
\label{G_3.int}
\begin{array}{l}
\displaystyle 2 \log \Gamma_{3}(z+1) = -\frac{z^3}{2}(\log z -
H_{3}) + \log G(z+1) \> + \\
\noalign{\vskip1.0pc}{\hskip1.5pc} \displaystyle

\Big(z^{2} \zeta^{\prime}(0) -2 z \zeta^{\prime}(-1) +
\zeta^{\prime}(-2) \Big) + 2 \int_{0}^{\infty} \frac{x^{2}
\arctan(x/z)}{e^{2\,\pi\,x}-1}\>dx
\end{array}
\end{equation}

\smallskip
\noindent that follows directly from Proposition 3.1 with $n=1$.
Expanding $arctan$ into the Taylor series

\[
\arctan(\frac{x}{z}) = \sum_{k=1}^{n} \frac{(-1)^{k-1}}{2k-1}
{(\frac{x}{z})}^{2k-1} + \frac{(-1)^n}{z^{2n-1}} \int_{0}^{x}
\frac{y^{2n}}{y^{2}+z^{2}} \> dy
\]

\noindent and making use of (\ref{int.bern}), we readily obtain

\[
\begin{array}{l}
\displaystyle \int_{0}^{\infty} \frac{x^{2}
\arctan(x/z)}{e^{2\,\pi\,x}-1}\>dx = -\sum_{k=1}^{n}
\frac{B_{2k+2}}{4\,(k+1)(2k-1) \, z^{2k-1}} \> + \\
\noalign{\vskip1.0pc}{\hskip7.0pc} \displaystyle

\frac{(-1)^n}{z^{2n-1}} \int_{0}^{\infty} \frac{t^{2}\>
dt}{e^{2\,\pi\,t}-1} \int_{0}^{t} \frac{y^{2n}}{y^{2}+z^{2}}\> dy.
\end{array}
\]

\noindent The conclusion of the lemma follows immediately by
substituting this into formula (\ref{G_3.int}) and using Lemma
5.1.

\end{proof}
\bigskip

\section{Implementation}

In this section we discuss a numeric computational scheme for the
Barnes function. Let us consider an integral
representation dated back to Alexejewsky \cite{Alexejewsky} and Barnes \cite{Barnes99}:

\[
\int_{0}^{z}\log  \Gamma
(x)d  x  = \frac{z(1-z)}{2}+\frac{z}{2}\log(2 \pi )+z \log \Gamma
(z)-\log  G(z+1),
\]
where $\Re(z) >  -1$. Employing integration by parts, we derive

\begin{equation}
\label{int_psi}
\log  G(z+1)=\frac{z(1-z)}{2}+\frac{z}{2}\log(2 \pi )+
\int_{0}^{z}x \, \psi (x)\, dx, \quad \Re(z)>-1
\end{equation}
The integral representation (\ref{int_psi}) provides an efficient
numeric procedure for arbitrary precision evaluation of the ${\rm
G}$-function. This representation demonstrates that the complexity
of computing $G(z)$ depends at most on the computational
complexity of the polygamma function. The later can be numerically
computed by the Spouge approximation \cite{Spouge} (a modification
of the Lanczos approximation). For the numerical integration in
(\ref{int_psi}) we shall use the Gaussian quadrature scheme.

\medskip
The restriction $\Re(z)>-1$ in (\ref{int_psi}) can be easily
removed by analyticity of the polygamma. For example, by resolving
the singularity of the integrand at the pole $x=-1$, we
continue \(\log G(z+1)\) to the wider area $\Re(z)>-2$:

\[
\log G(z+1)=\frac{z(1-z)}{2}+\frac{z}{2}\log(2 \pi )+\log  (z+1)  +\int _{0}^{z}\left(x \psi (x)-\frac{1}{x+1}\right)  dx
\]
Another method of the analytic continuation of the $G$-function is to
use the following well-known identity

\[
\psi (x)=\psi (-x)-\pi \cot(\pi x)-\frac{1}{x}
\]
which, upon substituting it into (\ref{int_psi}), yields

\begin{equation}
\label{kinkelin}
\log G(1-z)=\log G(1+z)-z  \log(2 \pi )+\int_{0}^{z}\pi   x \cot  (\pi x) \> dx.
\end{equation}
The identity (\ref{kinkelin}) holds everywhere in a complex plane
of $z$, except the real axes, where the integrand has simple
poles. Therefore, combining (\ref{int_psi}) with (\ref{kinkelin})
yields the following definition of the Barnes function

\begin{equation}
\label{contour}
\begin{array}{l}
\displaystyle
G(-n)=0, \quad n \in \Nset, \\
\noalign{\vskip1.0pc} \displaystyle
G(z)=(2 \pi)^{(z-1)/2} \exp \left(-\frac{
(z-1)(z-2)}{2}+{{\int}_{\gamma}}x \, \psi (x)\, dx\right),  \\
\noalign{\vskip1.0pc}{\hskip20.0pc} \displaystyle
\arg(z) \neq \pi,
\end{array}
\end{equation}
where the contour of integration \(\gamma \) is a line between 0
and $z-1$ that does not cross the negative real axis; for example,
\(\gamma \) could be the following path  \(\{0,i,i+z-1,z-1\}\).

It remains to define $G(z)$ for $z \in \Rset^{-}$. We do
this by using the reflexion formula (\ref{reflex}), which upon
periodicity of the Clausen function (\ref{clausen}) can be
rewritten as follows:

\begin{equation}
\label{G_clausen}
\begin{array}{l}
\displaystyle
G(-z)=(-1)^{{\lfloor z/2 \rfloor} -1}
G(z+2) {\Big| \frac{\sin(\pi z)}{\pi}\Big|}^{z+1} * \\
\noalign{\vskip1.0pc}{\hskip4.0pc} \displaystyle
\exp \Big(\frac{1}{2 \pi}{{\rm Cl}_2}\big(2 \pi (z-\lfloor z\rfloor)\big)
\Big), \quad z \in \Rset^{-}.
\end{array}
\end{equation}

\bigskip
\bigskip

\section{Glaisher's constant }

In this section we discuss numeric computational schemes for
Glaisher's constant:

\begin{equation}
\label{glaisher}
 A=\exp\Big(\frac{1}{12}-{{\zeta }^{\prime
}}(-1)\Big)
\end{equation}
In light of the functional equation for the zeta function

\[
\zeta (s)={2^s}{{\pi }^{s-1}}\Gamma (1-s)  \zeta
(1-s)\sin\Big(\frac{\pi s}{2}\Big)
\]
\smallskip
we can rewrite (\ref{glaisher}) as

\begin{equation}
\label{glaisher2}
A=\exp\Big(\frac{\gamma }{12}-\frac{{{\zeta
}^{\prime }}(2)}{2 {{\pi }^2}}\Big)  {{(2  \pi )}^{1/12}}.
\end{equation}
\smallskip
Exactly this formula is used in {\itshape Mathematica} V4.2 for
computing Glaisher's constant.

\medskip
Another computational scheme for Glaisher's constant is based on the Barnes $G$-function identity:

\begin{equation}
 \log A = \frac{1}{12} + \frac{\log 2}{36} - \frac{\log \pi
}{6} - \frac{2}{3} \> \log G\Big(\frac{1}{2}\Big).
\end{equation}
\noindent
Here $G(\frac{1}{2})$ can be computed via integral (\ref{int_psi}).

\medskip
    A more efficient algorithm is provided by the following identity

\begin{equation}
\label{A_series}
  \log A=\frac{\log 2}{12}+ \frac{1}{36}\sum _{
k=1}^{\infty } {\Big(\zeta (2 k+1)-1\Big)}
\Big(28+\frac{3}{1+k}-\frac{6}{2+k}\Big),
\end{equation}
which allows to approximate Glaisher's constant as

\[
\log A=\frac{\log 2}{12}+ \frac{1}{36}\sum _{
k=1}^{N}{\Big(\zeta (2 k+1)-1\Big)}
\Big(28+\frac{3}{1+k}-\frac{6}{2+k}\Big)+
O\Big(\frac{1}{{4^N}}\Big)
\]
This approximation is suitable for arbitrary
precision computation - it requires \(\lceil \frac{p}{2} \>
{{\log}_2}10\rceil \) terms to achieve $p$ decimal-digit accuracy.

\medskip
We close this section with a proof of formula (\ref{A_series}).
Recall that \cite{Bateman}

\begin{equation}
\label{zeta-1} \zeta (s) - 1 = \frac{1}{\Gamma (s)}
\int_{0}^{\infty} \frac{x^{s-1}} {e^x (e^{x}-1)} \> dx, \quad
\Re(s)>1.
\end{equation}

\smallskip
\noindent
Upon replacing $\zeta (2 k+1)-1$ in (\ref{A_series}) by its
integral representation (\ref{zeta-1}) and inverting the order of
summation and integration, we obtain

\begin{equation}
\label{A_series2} \log A=\frac{\log 2}{12}+ \frac{1}{36}
\int_{0}^{\infty} \frac{dx} {e^x (e^{x}-1)} \sum _{ k=1}^{\infty }
\frac{t^{2k}}{(2k)!} \> \Big(28+\frac{3}{1+k}-\frac{6}{2+k}\Big).
\end{equation}
Next, noting that
\[
\sum_{k=1}^{\infty } \frac{t^{2k}}{(2k)!} = \cosh (t) - 1
\]
and

\[
\sum_{k=1}^{\infty } \frac{t^{2k}}{(2k)! \> (k+s)} =
\frac{2}{t^{2s}} \int_{0}^{t} x^{2s-1} (\cosh (x) - 1) \> dx,
\]

\medskip
\noindent
we compute the inner sum

\[
\begin{array}{l}
\displaystyle
 \frac{t^4}{2}\sum _{ k=1}^{\infty }
\frac{t^{2k}}{(2k)!} \Big(28+\frac{3}{1+k}-\frac{6}{2+k}\Big) = -36 + 3 \,t^2 - 14 \,t^4 \> + \\
\noalign{\vskip1.0pc}{\hskip3.0pc} \displaystyle

  \left( 36 + 15 \,t^2 + 14 \,t^4 \right) \,
   \cosh(t) - \left( 36\, t + 3\, t^3 \right) \,
   \sinh(t)
\end{array}
\]

\smallskip
\noindent
In the next step, we substitute this back to (\ref{A_series2})
and, employed by (\ref{hurwitz_int}), evaluate the integral.
Unfortunately, integration cannot be done term by term. The way around
is to multiply the integrand by $t^\lambda, \,\lambda > 0$ and
consider the limit as $\lambda \rightarrow 0$.  This yields

\[
\begin{array}{l}
\displaystyle
 \int_{0}^{\infty} \Bigg(
 \frac{36\,t^{\lambda -4 }}{e^{t} - 1} +
  \frac{36\,t^{\lambda -4 }}
   {e^{2\,t}\,\left( e^{t} - 1 \right) } -
  \frac{72\,t^{\lambda -4 }}
   {e^t\,\left( e^{t} - 1 \right) } -
  \frac{36\,t^{\lambda -3 }}{e^{t} - 1} \>+\\
\noalign{\vskip1.0pc} {\hskip2.5pc}\displaystyle

  \frac{36\,t^{\lambda -3 }}
   {e^{2\,t}\,\left( e^{t} - 1 \right) } +
  \frac{15\,t^{\lambda -2 }}{e^{t} - 1} +
  \frac{15\,t^{\lambda -2 }}
   {e^{2\,t}\,\left( e^{t} - 1 \right) } +
  \frac{6\,t^{\lambda -2 }}
   {e^t\,\left( e^{t} - 1 \right) } -
  \frac{3\,t^{\lambda -1 }}{e^{t} - 1} \>+ \\
\noalign{\vskip1.0pc}{\hskip2.5pc} \displaystyle

  \frac{3\,t^{\lambda -1 }}
   {e^{2\,t}\,\left( e^{t} - 1 \right) } +
  \frac{14\,t^{\lambda }}{e^{t} - 1} +
  \frac{14\,t^{\lambda }}
   {e^{2\,t}\,\left( e^{t} - 1 \right) } -
  \frac{28\,t^{\lambda }}
   {e^t\,\left( e^{t} - 1 \right) } \Bigg) \> dt =\\
\noalign{\vskip1.0pc} \displaystyle
   36\,\Gamma(\lambda -3 ) -
  9\,2^{5 - \lambda }\,
   \Gamma(\lambda -3 ) -
  36\,\Gamma(\lambda -2 ) -
  9\>2^{4 - \lambda }\,
   \Gamma(\lambda -2 ) \>- \\
\noalign{\vskip1.0pc} \displaystyle
  21\,\Gamma(\lambda -1 ) -
  15\,2^{1 - \lambda }\,
   \Gamma(\lambda -1 ) -
  3\,\Gamma(\lambda ) -
  \frac{3\,\Gamma(\lambda )}
   {2^{\lambda }} +
  14\,\Gamma(\lambda +1 ) \>-\\
\noalign{\vskip1.0pc} \displaystyle
  \frac{7\,\Gamma(\lambda +1 )}
   {2^{\lambda }} +
  36\,\Gamma(\lambda -1 )\,
   \zeta(\lambda -1 ).
\end{array}
\]

\smallskip
\noindent
Computing the limit as $\lambda \rightarrow 0$, we find that the right hand-side of the above expression evaluates to

\[3 - 3 \log 2 - 36 \,\zeta^{\prime}(-1).
\]

\medskip
\noindent
This completes the proof of (\ref{A_series}).

\end{document}